\documentclass[10pt,reqno]{amsart}
\usepackage[top=30truemm,bottom=30truemm,left=25truemm,right=25truemm]{geometry}

\usepackage{amssymb, amsmath, amsthm, amsfonts, latexsym, mathtools, amsfonts}
\usepackage{graphicx}
\usepackage{xcolor}
\usepackage{ytableau}
\usepackage{tikz, comment}
\usepackage{amscd}
\usepackage{dsfont, enumitem}
\usepackage{arydshln}
\usepackage[normalem]{ulem}
\usepackage{mathtools}
\usepackage{thm-restate}
\usepackage{booktabs}
\usepackage{tikz-cd}

\makeatletter
\@ifundefined{newcounteralias}{}{%
  \renewcommand\thmt@autorefsetup{%
    \@xa\def\csname\thmt@envname autorefname\@xa\endcsname\@xa{\thmt@thmname}%
  }%
}
\makeatother

\usepackage{hyperref}
\usepackage[nameinlink]{cleveref}
\hypersetup{
setpagesize=false,
bookmarksnumbered=true,%
bookmarksopen=true,%
colorlinks=true,%
linkcolor=orange,
citecolor=blue,
}

\crefname{equation}{}{}

\newtheorem{thm}{Theorem}[section]
\crefname{thm}{Theorem}{Theorems}

\crefname{prop}{Proposition}{Propositions}
\newtheorem{lem}[thm]{Lemma}
\crefname{lem}{Lemma}{Lemmas}
\newtheorem{cor}[thm]{Corollary}
\crefname{cor}{Corollary}{Corollarys}
\newtheorem{conj}[thm]{Conjecture}
\crefname{conj}{Conjecture}{Conjectures}

\theoremstyle{definition}
\newtheorem*{definition}{Definition}

\newtheorem*{remark}{Remark}

\crefname{section}{Section}{Sections}

\newcommand{\til}{\widetilde}

\newcommand{\Img}{\operatorname{Im}}

\newcommand{\ceq}{\coloneq}

\newcommand{\SL}{\operatorname{SL}}

\newcommand{\wt}{\operatorname{wt}}

\newcommand{\CC}{\mathbb{C}}
\newcommand{\RR}{\mathbb{R}}
\newcommand{\QQ}{\mathbb{Q}}
\newcommand{\ZZ}{\mathbb{Z}}

\newcommand{\HH}{\mathbb{H}}
\newcommand{\PP}{\mathbb{P}}

\newcommand{\E}{\mathcal{E}}
\newcommand{\Ead}{\mathcal{E}^{\mathrm{adm}}}
\newcommand{\EY}{\mathcal{E}^{Y}}
\newcommand{\tEY}{\mathcal{\til{E}}^{Y}}

\newcommand{\W}{\mathcal{W}}
\newcommand{\Q}{\mathcal{Q}}

\font\sevency=wncyr7
\newcommand{\sha}{\hbox{\sevency X}}

\numberwithin{equation}{section}

\begin{document}

\title[On the Space Defined by the Modular Transformation Law for Multiple Eisenstein Series]{On the Space Defined by the Modular Transformation Law\\for Multiple Eisenstein Series}
\author{Katsumi Kina}
\address{Graduate School of Mathematics, Kyushu University, Motooka 744, Nishi-ku, Fukuoka 819-0395 Japan}
\email{kkina.math@gmail.com}

\begin{abstract}
Motivated by a modular transformation law for multiple Eisenstein series established by Yu, we introduce spaces of holomorphic functions satisfying this law together with an additional filtration condition. We show that these spaces contain the spaces spanned by admissible multiple Eisenstein series, are finite-dimensional in each weight, and share several structural properties with them. We also obtain an explicit formula for the generating function of their dimensions. However, the inclusions are known to be strict for certain weights. Thus, the defining conditions do not suffice to characterize the spaces spanned by admissible multiple Eisenstein series.
\end{abstract}

\maketitle


\section{Introduction and Definitions}

Multiple Eisenstein series are holomorphic functions on the upper half-plane $\HH$, introduced by Gangl, Kaneko, and Zagier \cite{GKZ06} and defined by the following iterated sums. For an integer $r>0$ and integers $k_1,\dots,k_r\ge2$, we define the multiple Eisenstein series by
\begin{equation}\label{eq:MES-def}
G(k_1,\dots,k_r)(\tau)\ceq
\lim_{M\to\infty}\lim_{N\to\infty}\sum_{\substack{w_s\in\ZZ_M\tau+\ZZ_N\\ 0\prec w_1\prec\dots\prec w_r}}\frac{1}{w_1^{k_1}w_2^{k_2}\cdots w_r^{k_r}},
\end{equation}
Here, $\ZZ_N=\{n\in\ZZ\mid |n|<N\}$ and the order $\prec$ on the lattice points in $\ZZ\tau+\ZZ$ is defined by
\begin{align}
m_1\tau+n_1\prec m_2\tau+n_2 &\overset{\mathrm{def}}{\Longleftrightarrow}
\left\{\begin{array}{c}
m_1<m_2
\\
\text{or}
\\
m_1=m_2 \quad \text{and}\quad n_1<n_2
\end{array}\right. .
\end{align}
We call a tuple $(k_1,\dots,k_r)$ an index, its weight $k_1+\cdots+k_r$, and its depth $r$.
An index is said to be admissible if $k_r\ge 3$.

We are interested in the structure of the space spanned by multiple Eisenstein series. For each positive integer $k$, we define
$$\Ead_k \ceq \langle G(\mathds{k}) \mid \mathds{k}\text{ is an admissible index of weight }k\rangle_{\RR}$$
and
$$\E_k \ceq \langle G(\mathds{k}) \mid \mathds{k}\text{ is an index of weight }k \rangle_{\RR}.$$
We also set $\Ead_0=\E_0=\RR$.

Although spaces of multiple Eisenstein series are usually considered over $\QQ$, we work over $\RR$ throughout this paper. By \cite[Theorem B (ii)]{BK26}, the corresponding spaces over $\RR$ and $\CC$ are obtained from those over $\QQ$ by extension of scalars, so no new linear relations arise upon extending the coefficient field.
Both $\sum_k\E_k$ and $\sum_k\Ead_k$ are rings under multiplication, with the product given by the harmonic product. In particular,
$$\E_r\cdot\E_s\subseteq\E_{r+s}
\qquad\text{and}\qquad
\Ead_r\cdot\Ead_s\subseteq\Ead_{r+s}.$$
Moreover, the following result is known.

\begin{thm}[{\cite[Theorem A (iii), Theorem B (i)]{BK26}}]\label{thm:weight-graded}
We have
$$\sum_{k\geq 0}\E_k = \bigoplus_{k\geq 0}\E_k
= \bigoplus_{a,k\geq 0}G(2)^a\Ead_k.$$
This shows that $\sum_{k\geq 0}\E_k$ is a graded ring
and that $G(2)$ is transcendental over $\sum_{k\geq 0}\Ead_k$.
\end{thm}

In this article, we focus particularly on the following two results of \cite{Yu26}. For a function $f$ and an integer $k$, we define
$$f^{\flat_k}(\tau) \ceq (-1)^k \overline{f(-\overline{\tau})}
\qquad\text{and}\qquad
\sigma_k(f)(\tau) \ceq f(\tau) - f^{\flat_k}(\tau).$$
It is easy to see that
\begin{equation}\label{eq:inv-flat}
\sigma_k(f)^{\flat_k} = - \sigma_k(f).
\end{equation}
Moreover, if $f(\tau)=\sum_{n}a_nq^n$, where $q=e^{2\pi i\tau}$, then
\begin{equation}\label{eq:flat-Fourier}
f(\tau)^{\flat_k} = (-1)^k \sum_{n}\overline{a_n}q^n.
\end{equation}
Throughout this paper, for each $u\in\QQ$, we choose a matrix $\gamma_u\in\mathrm{SL}_2(\ZZ)$ whose bottom row is $(c,d)$ with $c\ge1$ and $d/c=u$. Although $\gamma_u$ is not uniquely determined, our subsequent arguments are independent of this choice.

\begin{thm}[{\cite[Theorem 1.5]{Yu26}}]\label{thm:modular-trans}
For $f\in\Ead_k$, we have
\begin{equation}\label{eq:Yu-trans}
f|_k\gamma_u
= (-1)^kf + \sum_{\substack{w\in\QQ\\w\leq u}} \sigma_k(f)|_k\gamma_{w}.\tag{Yu}
\end{equation}
The series is absolutely convergent on $\mathbb{H}$, and each summand is independent of the choice of $\gamma_{w}$.
\end{thm}

\begin{thm}[{\cite[Lemma 4.6]{Yu26}}]\label{thm:weight-descent}
For $f\in\E_k$, we have
$$\sigma_k(f) \in \E_{0}+\E_{1}+\cdots+\E_{k-1}.$$
\end{thm}

We refer to the transformation law \cref{eq:Yu-trans} in \cref{thm:modular-trans} as Yu transformation law of weight $k$. Motivated by Yu's results, we introduce the following spaces.

\begin{definition}
For each integer $m \geq 0$, set
$$\W_{-m-1} \ceq 0
\qquad\text{and}\qquad
\W_{m} \ceq \sum_{\substack{a,j\geq 0\\2a+j\leq m}} G(2)^a\EY_j.$$
For each integer $k$, we inductively define $\EY_k$ to be the $\RR$-vector space of functions $f$ satisfying the following conditions:
\begin{enumerate}[label=\textbf{(C\arabic*)}]
\setlength{\itemsep}{2pt}
    \item $f$ is holomorphic on $\HH$, is $1$-periodic, and is holomorphic at the cusp $i\infty$.
    \item $f$ satisfies Yu transformation law \cref{eq:Yu-trans} of weight $k$, with the series converging absolutely.
    \item $\sigma_k(f)\in\W_{k-1}$.
\end{enumerate}
We call $\EY_k$ the space of Yu transformation forms of weight $k$. An element of $\EY_k$ is called a Yu transformation form of weight $k$.
\end{definition}

Let
$$T \ceq \begin{pmatrix}
1 & 1\\
0 & 1
\end{pmatrix}
\qquad\text{and}\qquad
S \ceq \begin{pmatrix}
0 & -1\\
1 & 0
\end{pmatrix}.$$

\begin{remark}
The independence of Yu transformation law in \textbf{(C2)} from the choice of $\gamma_w$ follows from the $1$-periodicity of $f$ in \textbf{(C1)}. Indeed, $\sigma_k(f)$ is also $1$-periodic, and any two choices $\gamma_w$ and $\gamma'_w$ differ by left multiplication by a power of $T$. Hence, each summand in \cref{eq:Yu-trans} is independent of the choice of $\gamma_w$.
\end{remark}

The following result is an immediate consequence of \cref{thm:weight-graded,thm:modular-trans,thm:weight-descent}.
\begin{thm}\label{thm:inclusion}
For each $k\ge 0$, we have $\Ead_k\subseteq\EY_k$.
\end{thm}

In \cref{sec:prop-EY}, we study several properties of $\EY_k$. In \cref{sec:upper_dim}, we establish a dimension formula for $\EY_k$.
It is known the inclusion $\Ead_k\subseteq\EY_k$ is strict for some weights. This raises the natural question of whether additional conditions can be imposed on the definition of $\EY_k$ to obtain a space coinciding with $\Ead_k$. At present, however, we do not know of any such conditions.

\section*{Use of AI}
The author used ChatGPT 6 as a supplementary tool in developing some of the mathematical statements and proofs presented in this paper. Several ideas, computations, and arguments were developed through discussions with the model. All mathematical statements and final proofs were independently verified by the author, who takes full responsibility for the content of this paper.

\section*{Acknowledgments}
The author expresses his gratitude to his academic advisor Professor Masanobu Kaneko and Professor Hiroyuki Ochiai.
The author would also like to thank Professor Henrik Bachmann for valuable advice on this research.

\section{Properties of $\EY_k$}\label{sec:prop-EY}

\subsection{Basic Properties of $\EY_k$}

Throughout this paper, \(M_k\) and \(S_k\) denote the real vector spaces of holomorphic modular forms and cusp forms of weight \(k\) for \(\SL_2(\ZZ)\), respectively, whose Fourier coefficients are real.
It is well known that $M_k$ is spanned over $\RR$ by the monomials $G(4)^aG(6)^b$ with $4a+6b=k$.

\begin{lem}\label{lem:MES-modular}
For $f\in\EY_k$, $f\in M_k$ if and only if $\sigma_k(f)=0$.
\end{lem}

\begin{proof}
If $f\in M_k$, then $\sigma_k(f)=0$ since $f$ has real Fourier coefficients.
Conversely, suppose that $\sigma_k(f)=0$.
If $k$ is even, then \textbf{(C2)} implies that $f\in M_k$.
If $k$ is odd, applying $S$ twice to $f|_kS=-f$ yields $f=-f$, so $f=0\in M_k$.
\end{proof}

The following theorem motivates the inclusion of condition \textbf{(C3)}, in addition to the transformation law, in the definition of $\EY_k$.

\begin{thm}\label{thm:finite}
For each $k\geq 0$, we have $\dim\EY_k<\infty$.
\end{thm}

\begin{proof}
By \cref{lem:MES-modular}, the kernel of $\sigma_k:\EY_k\to\W_{k-1}$ is $M_k$. Hence, we obtain
\begin{equation}\label{eq:rough-bound1}
\dim\EY_k\le \dim M_k + \dim\W_{k-1}.
\end{equation}
Since $\W_{k-1}$ is finite-dimensional by induction on $k$, it follows that $\EY_k$ is also finite-dimensional.
\end{proof}

\begin{lem}\label{lem:limit-f}
For $f\in \EY_k$, we have
$$\lim_{u\to-\infty}(f|_k\gamma_u)(\tau) = (-1)^kf(\tau).$$
\end{lem}

\begin{proof}
By \textbf{(C2)}, we have
$$(f|_k\gamma_u)(\tau) = (-1)^kf(\tau) + \sum_{w\leq u}(\sigma_k(f)|_k\gamma_w)(\tau).$$
Since the series converges absolutely, its tail tends to zero as $u\to-\infty$, yielding the desired result.
\end{proof}

\begin{lem}\label{lem:flat-mod-inv}
For any integers $k$ and $m$, and $f\in\W_m$, we have
$$f^{\flat_k}\equiv (-1)^{k-m}f\pmod{\W_{m-1}}.$$
\end{lem}

\begin{proof}
Note that $G(2)^{\flat_2}=G(2)$. For integers $a,j\ge0$ with $2a+j\le m$ and $h\in\EY_j$, we have
$$(G(2)^ah)^{\flat_k} = (-1)^{k-j}G(2)^ah^{\flat_j} = (-1)^{k-j}G(2)^a(h-\sigma_j(h)).$$
Since $\sigma_j(h)\in\W_{j-1}$, we have $G(2)^a\sigma_j(h)\in\W_{m-1}$. Moreover, if $2a+j=m$, then $(-1)^{k-j}=(-1)^{k-m}$, while if $2a+j<m$, then $G(2)^ah\in\W_{m-1}$. This proves the assertion.
\end{proof}

\begin{lem}\label{lem:strong-cond}
For any $f\in\EY_k$, we have $\sigma_k(f)\in\W_{k-3}$. In other words, condition \textbf{(C3)} is equivalent to the stronger condition $\sigma_k(f)\in\W_{k-3}$.
\end{lem}

\begin{proof}
Since $\sigma_k(f)\in\W_{k-1}$, we may write
$$\sigma_k(f)=\sum_{m=0}^{k-1}H_m, \qquad
H_m=\sum_{\substack{a,j\ge 0\\2a+j=m}}G(2)^a h_{a,j}, \qquad
h_{a,j}\in\EY_j.$$
By the modular transformation formula for $G(2)$,
$$\lim_{u\to-\infty}(G(2)|_2\gamma_u)(\tau)
=\lim_{u\to-\infty}\left(G(2)(\tau)-\frac{\pi i}{\tau+u}\right) = G(2)(\tau).$$
Together with \cref{lem:limit-f}, this yields
$$\lim_{u\to-\infty}(H_m|_m\gamma_u)(\tau) = (-1)^mH_m(\tau).$$
Consequently,
$$\lim_{u\to-\infty}(c\tau+d)(\sigma_k(f)|_k\gamma_u)(\tau)
= (-1)^{k-1}H_{k-1}(\tau),$$
where $u=d/c$, $c\ge1$, and $\gcd(c,d)=1$. Suppose that $H_{k-1}\neq0$, and choose $\tau_0\in\HH$ such that $H_{k-1}(\tau_0)\neq0$. Then
$$(\sigma_k(f)|_k\gamma_u)(\tau_0)\sim\frac{C_0}{c\tau_0+d}, \qquad C_0\ceq(-1)^{k-1}H_{k-1}(\tau_0)\neq0,$$
as $u\to-\infty$. This contradicts the absolute convergence of
$$\sum_{u\le N}(\sigma_k(f)|_k\gamma_u)(\tau_0).$$
Therefore, $H_{k-1}=0$, so $\sigma_k(f)\in\W_{k-2}$. Although the vanishing of $H_{k-2}$ can be proved similarly, we give an alternative proof here. By \cref{lem:flat-mod-inv,eq:inv-flat}, we have
$$-\sigma_k(f)=\sigma_k(f)^{\flat_k} \equiv \sigma_k(f) \mod\W_{k-3}.$$
Hence $\sigma_k(f)\in\W_{k-3}$, as desired.
\end{proof}

\Cref{lem:strong-cond} gives a sharper upper bound on $\dim\EY_k$ than \cref{eq:rough-bound1}. Indeed, since the kernel of $\sigma_k:\EY_k\to\W_{k-3}$ is $M_k$, we obtain
\begin{equation}\label{eq:rough-bound2}
\dim\EY_k\le \dim M_k+\dim\W_{k-3}.
\end{equation}
Combining this bound with \cref{thm:inclusion}, we obtain
$$\EY_k=0\quad(k<0),\qquad \EY_0=\RR,\qquad \EY_1 = \EY_2 = 0,$$
$$\EY_3=\RR\cdot G(3),\qquad \EY_4=\RR\cdot G(4),\qquad\EY_5=\RR\cdot G(5)+\RR\cdot G(2,3).$$

\subsection{Properties of $\Ead_k$ Extending to $\EY_k$}\label{subsec:ext-prop}

We now show that several known properties of $\Ead_k$ also hold for $\EY_k$ (see \cite{BKM26,BK26}). Our results are summarized as follows.

\begin{enumerate}[leftmargin=30mm]
\setlength{\itemsep}{6pt}
\item[(\cref{thm:cusp-limit})] For $f\in\EY_k$, $(f|_k\gamma_u)(i\infty)=(-1)^kf(i\infty)$. Furthermore, $f(i\infty)\in\RR$ by \cref{lem:const-real}.
\item[(\cref{thm:product})] $\EY_r\cdot\EY_s\subseteq\EY_{r+s}$.
\item[(\cref{thm:graded})] $\EY\ceq \sum_{k}\EY_k=\bigoplus_{k}\EY_k$.
\item[(\cref{thm:G2-indep})] $\EY[G(2)]\cong \EY[X]$.
\item[(\cref{thm:serre-der})] $\partial_k(\EY_k)\subseteq\EY_{k+2}$, where $\partial_k$ denotes the Serre derivative.
\end{enumerate}

\begin{lem}\label{lem:adm-ord}
For $f\in\EY_k$, the following hold:
\begin{enumerate}
\item As $u\to\infty$,
$$(f|_k\gamma_{-u})(\tau) = (-1)^kf(\tau) + O(u^{-2}),$$
locally uniformly in $\tau\in\HH$.
\item As $\operatorname{Im}(\tau)\to\infty$,
$$(f|_k\gamma_u)(\tau) = (-1)^kf(\tau) + O(\operatorname{Im}(\tau)^{-2}),$$
uniformly in $u\in\QQ$ and $\operatorname{Re}(\tau)\in\RR$.
\end{enumerate}
\end{lem}

\begin{proof}
We argue by induction on $k$, showing simultaneously that $(f|_k\gamma_u)(\tau)$ is uniformly bounded in $u\in\QQ$,
locally uniformly in $\tau\in\HH$, and uniformly for $\operatorname{Im}(\tau)\ge1$. The case $k\le2$ is immediate. Let $k\ge3$ and assume the assertions for all smaller weights. For $h\in\W_{k-3}$, we may write
$$h=\sum_{2a+j\le k-3}G(2)^a h_{a,j}, \qquad h_{a,j}\in\EY_j.$$
Using the transformation formula for \(G(2)\), the induction hypothesis, and \(2a+j\le k-3\), we obtain
\begin{equation}\label{eq:ord-3}
(h|_k\gamma_{d/c})(\tau) = O(|c\tau+d|^{-3}),
\end{equation}
uniformly in $d/c\in\QQ$, locally uniformly in $\tau\in\HH$, and uniformly for $\operatorname{Im}(\tau)\ge1$. Now \textbf{(C2)} gives
$$(f|_k\gamma_u)(\tau)-(-1)^kf(\tau) = \sum_{w\le u}(\sigma_k(f)|_k\gamma_w)(\tau).$$
Since $\sigma_k(f)\in\W_{k-3}$ by \cref{lem:strong-cond}, the estimate \cref{eq:ord-3}, together with the elementary estimates
$$\sum_{\substack{c\ge1,\,d\in\ZZ\\d/c\le-u}}|c\tau+d|^{-3}=O(u^{-2})$$
locally uniformly in $\tau$, and
$$\sum_{c\ge1,\,d\in\ZZ}|c\tau+d|^{-3} = O(\operatorname{Im}(\tau)^{-2})$$
uniformly in $\operatorname{Re}(\tau)$ as $\operatorname{Im}(\tau)\to\infty$, proves (1) and (2). Finally, since $f$ is bounded on $\operatorname{Im}(\tau)\ge1$ by \textbf{(C1)}, the same estimates imply the required uniform boundedness of $(f|_k\gamma_w)(\tau)$. This completes the induction.
\end{proof}

\begin{thm}\label{thm:cusp-limit}
Let $f\in\EY_k$ and $\alpha\ceq\lim_{\tau\to i\infty}f(\tau)$. Then, we have
$$\lim_{\tau\to i\infty}(f|_k\gamma_{u})(\tau)=(-1)^k\alpha.$$
\end{thm}

\begin{proof}
This follows immediately from \cref{lem:adm-ord}(2).
\end{proof}

\begin{thm}\label{thm:product}
For any integers $r$ and $s$, we have $\EY_r\cdot\EY_s\subseteq\EY_{r+s}$.
\end{thm}

\begin{proof}
For $f\in\EY_r$, from \textbf{(C2)}, we have
\begin{equation}\label{eq:flat-trans}
f^{\flat_r}|_r\gamma_u = (f - \sigma_r(f))|_r\gamma_u
= (-1)^rf + \sum_{w<u} \sigma_r(f)|_r\gamma_{w}.
\end{equation}
Thus, for $f\in\EY_r$ and $g\in \EY_s$, we have
\begin{align*}
\sigma_{r+s}(fg) &= fg - f^{\flat_r}g^{\flat_s}
= (f^{\flat_r}+\sigma_r(f))g - f^{\flat_r}(g-\sigma_s(g))
= \sigma_r(f)g + f^{\flat_r}\sigma_s(g)
\end{align*}
Applying $|_{r+s}\gamma_u$ to the above equation and substituting \cref{eq:flat-trans} and \textbf{(C2)}, we obtain
\begin{align*}
(\sigma_{r+s}(fg))|_{r+s}\gamma_u
&= \sigma_r(f)|_r\gamma_u\left((-1)^sg + \sum_{w\leq u} \sigma_s(g)|_s\gamma_w\right)
+ \left((-1)^rf + \sum_{w< u} \sigma_r(f)|_r\gamma_w\right)\sigma_s(g)|_s\gamma_u
\end{align*}
Since the series in \textbf{(C2)} for $f$ and $g$ converge absolutely, all the following sums converge absolutely. Summing over $u\le N$, we obtain
\begin{align*}
\sum_{u\leq N}(\sigma_{r+s}(fg))|_{r+s}\gamma_u
&= (-1)^sg\sum_{u\leq N}\sigma_r(f)|_r\gamma_u + \sum_{w\leq u\leq N} \sigma_r(f)|_r\gamma_u\sigma_s(g)|_s\gamma_w
\\
&\qquad\qquad + (-1)^rf\sum_{u\leq N}\sigma_s(g)|_s\gamma_u + \sum_{w< u\leq N} \sigma_r(f)|_r\gamma_w\sigma_s(g)|_s\gamma_u.
\end{align*}
From the identity $\sum_{w\leq u\leq N}+\sum_{u<w\leq N}=\sum_{u,w\leq N}$, we obtain
\begin{align*}
\sum_{u\leq N}(\sigma_{r+s}(fg))|_{r+s}\gamma_u
&= (-1)^sg(f|_r\gamma_N-(-1)^{r}f) + (-1)^rf(g|_s\gamma_N-(-1)^{s}g)
\\
&\qquad\qquad + (f|_r\gamma_N-(-1)^{r}f)(g|_s\gamma_N-(-1)^{s}g)
\\
&= (fg)|_{r+s}\gamma_N - (-1)^{r+s}fg.
\end{align*}
Therefore, $fg$ satisfies \textbf{(C2)}. Using the identity
\begin{equation}\label{eq:prod-sigma}
\sigma_{r+s}(fg) = f\sigma_s(g)+\sigma_r(f)g-\sigma_r(f)\sigma_s(g)
\end{equation}
and induction, we conclude that $fg$ satisfies \textbf{(C3)}. This completes the proof.
\end{proof}

\begin{thm}\label{thm:graded}
$\EY\ceq\sum_k\EY_k$ is a graded ring.
\end{thm}

\begin{proof}
By \cref{thm:product}, it remains to show that the sum is direct. Let $k_1<\cdots<k_r$ be integers and suppose that
$$f_{k_1}+\cdots+f_{k_r}=0,
\qquad 0\neq f_{k_i}\in\EY_{k_i}.$$
Applying $|_{k_r}\gamma_u$ and letting $u\to-\infty$, we obtain $(-1)^{k_r}f_{k_r}=0$ by \cref{lem:limit-f}. Hence, the sum is direct.
\end{proof}

\begin{thm}\label{thm:G2-indep}
$G(2)$ is transcendental over $\EY$.
\end{thm}

\begin{proof}
We first prove the statement in the homogeneous case. Suppose that
$$P(X)\ceq\sum_{\substack{a,j\ge0\\2a+j=k}}f_jX^a\in\EY[X], \qquad f_j\in\EY_j,$$
satisfies $P(G(2))=0$.
By \cref{lem:adm-ord}(1) and the transformation formula for $G(2)$, we have
\begin{align*}
P(G(2))|_k\gamma_{-u}
&=(-1)^k\sum_{\substack{a,j\ge0\\2a+j=k}}(f_j+O(u^{-2}))\left(G(2)-\frac{\pi i}{\tau-u}\right)^a\\
&=(-1)^k\left(P(G(2))-\frac{\pi i}{\tau-u}P'(G(2)) + O(u^{-2})\right).
\end{align*}
Since $P(G(2))=0$, multiplying by $\tau-u$ and letting $u\to\infty$ yields $P'(G(2))=0$. Repeating the same argument for the successive derivatives of $P$, which are also homogeneous, we obtain
$$P^{(n)}(G(2))=0 \qquad(n\ge0).$$
This implies that $P(X)=0$.
For the general case, suppose that $P(G(2))=0$ and decompose $P(X)\in\EY[X]$ into its homogeneous components:
$$P(X)=\sum_{m=0}^k P_m(X),$$
where $P_m(X)$ is homogeneous of weight $m$. Then, we have
$$0 = P(G(2))|_k\gamma_u = \sum_{m=0}^k(c\tau+d)^{-k+m}\bigl(P_m(G(2))|_m\gamma_u\bigr)(\tau),$$
where $u=d/c$, $c\ge1$, and $\gcd(c,d)=1$. Since
$$\lim_{u\to-\infty}\bigl(P_m(G(2))|_m\gamma_u\bigr)(\tau) = (-1)^mP_m(G(2))(\tau),$$
letting $u\to-\infty$ yields $P_k(G(2))=0$. By the homogeneous case, we have $P_k(X)=0$. Repeating this argument for the remaining homogeneous components, we conclude that $P(X)=0$.
\end{proof}

For even $k\geq 2$, let $E_k\ceq G(k)/\zeta(k)$. Moreover, denote $q=e^{2\pi i\tau}$,
$$D=\frac{1}{2\pi i}\frac{d}{d\tau} = q\frac{d}{dq}
\ ,\quad \text{and}\quad
\partial_k = D-\frac{k}{12}E_2.$$

\begin{thm}\label{thm:serre-der}
We have $\partial_k(\EY_k)\subseteq\EY_{k+2}$.
\end{thm}

\begin{proof}
Let $f\in\EY_k$. We show that $\partial_k f$ satisfies conditions \textbf{(C1)}--\textbf{(C3)}.
Condition \textbf{(C1)} follows from the fact that $f$ is holomorphic, $1$-periodic, and holomorphic at $i\infty$, and hence has a Fourier expansion with no negative powers of $q$.
To verify \textbf{(C2)}, we apply the Serre derivative to both sides of Yu transformation law \cref{eq:Yu-trans}. Using the identities
$$\partial_k(f|_k\gamma)=(\partial_k f)|_{k+2}\gamma
\quad\text{and}\quad
\partial_k(\sigma_k(f))=\sigma_{k+2}(\partial_k f),$$
we obtain
\begin{align*}
(\partial_k f)|_{k+2}\gamma_u = \partial_k(f|_k\gamma_u)
= (-1)^k\partial_k f + \sum_{w\le u}\sigma_{k+2}(\partial_k f)|_{k+2}\gamma_w.
\end{align*}
Here, the interchange of differentiation and summation is justified by the locally uniform convergence of the series in \cref{eq:Yu-trans}, which follows from \cref{eq:ord-3}. Finally, we verify \textbf{(C3)} by induction on $k$. Suppose that the assertion holds for all weights less than $k$. By \cref{lem:strong-cond}, we have $\sigma_k(f) \in \W_{k-3}$. Using the Ramanujan identity $DE_2=(E_2^2-E_4)/12$, for $E_2^a h$ with $h\in\EY_j$ and $2a+j\le k-3$, we obtain
$$\partial_k(E_2^ah) = E_2^a\partial_jh + \frac{a+j-k}{12}E_2^{a+1}h - \frac a{12}E_2^{a-1}E_4h,$$
where the last term is understood to be zero if $a=0$. Since $j<k$, the induction hypothesis gives $\partial_jh\in\EY_{j+2}$. Moreover, $E_4h\in\EY_{j+4}$ by \cref{thm:product}. Thus, each term on the right belongs to $\W_{k-1}$. Hence,
\[
\partial_k(\W_{k-3})\subseteq\W_{k-1},
\]
and therefore
\[
\sigma_{k+2}(\partial_kf)
=\partial_k(\sigma_k(f))
\in\W_{k-1}.
\]
This proves \textbf{(C3)} and completes the proof.
\end{proof}

\section{The Dimension Formula for $\EY_k$}\label{sec:upper_dim}

Let $\W_{k-3}^-\ceq\{h\in\W_{k-3}\mid h^{\flat_k}=-h\}$.
By \cref{lem:strong-cond} and the identity $\sigma_k(f)^{\flat_k}=-\sigma_k(f)$, the $\RR$-linear map $\sigma_k:\EY_k\to\W_{k-3}^-$ is well-defined. Since its kernel is $M_k$, we obtain
$$\dim\EY_k=\dim M_k+\dim\sigma_k(\EY_k).$$
To determine the dimension of $\EY_k$, we study the image
of $\sigma_k$ using the following Poincar\'e series.
For a $1$-periodic function $f$, define
$$P_k(f) \ceq f+\sum_{u\in\QQ}f|_k\gamma_u,$$
whenever the series converges absolutely.

\begin{lem}\label{lem:poincare}
For $f\in\EY_k$, the series defining $P_k(\sigma_k(f))$ converges absolutely. Moreover, we have
$$P_k(\sigma_k(f))=(1-(-1)^k)f.$$
\end{lem}

\begin{proof}
The absolute convergence follows from \cref{eq:ord-3}. By \textbf{(C2)}, we have
$$f|_kS=(-1)^kf+\sum_{w\le0}\sigma_k(f)|_k\gamma_w.$$
Applying $S$ twice, we obtain
\begin{align*}
(-1)^kf = f|_kS^2 &= (-1)^kf|_kS + \sum_{w\leq 0}\sigma_k(f)|_k\gamma_wS
\\
&= f + (-1)^k\sum_{w\leq 0}\sigma_k(f)|_k\gamma_w + \sum_{w\leq 0}\sigma_k(f)|_k\gamma_wS
\\
&= f + (-1)^k\sigma_k(f) + (-1)^k\sum_{w\leq 0}\sigma_k(f)|_k\gamma_w + (-1)^k \sum_{w > 0}\sigma_k(f)|_k\gamma_w.
\end{align*}
Hence, we obtain
$$((-1)^k-1)f = (-1)^k P_k(\sigma_k(f)).$$
Multiplying both sides by $(-1)^k$ gives the desired identity.
\end{proof}

\subsection{The odd weight case}

\begin{thm}\label{thm:odd-dim}
For every odd integer $k$, the map
$$\sigma_k:\EY_k\longrightarrow\W_{k-3}^-$$
is an isomorphism of $\RR$-vector spaces, whose inverse is given by $\frac12P_k$. In particular,
$$\dim\EY_k = \dim\sigma_k(\EY_k) = \dim\W_{k-3}^-.$$
\end{thm}

\begin{proof}
Since $\ker\sigma_k=M_k=0$ by \cref{lem:MES-modular}, the map $\sigma_k$ is injective. To prove surjectivity, let $h\in\W_{k-3}^-$ and set
$$F \ceq \frac12P_k(h) = \frac12\left(h+\sum_{w\in\QQ}h|_k\gamma_w\right).$$
By \cref{eq:ord-3}, the series converges absolutely and locally uniformly on $\HH$, and $F$ has a finite limit at $i\infty$. Moreover, for each $w\in\QQ$, we have $\gamma_wT=T^n\gamma_{w+1}$ for some $n\in\ZZ$. Hence, $F$ satisfies \textbf{(C1)}. Since $k$ is odd, we have
$$(h|_k\gamma_w)^{\flat_k} = -h^{\flat_k}|_k\gamma_{-w}.$$
Thus, using $h^{\flat_k}=-h$, we obtain
$$P_k(h)^{\flat_k}
= h^{\flat_k} + \sum_{w\in\QQ}(h|_k\gamma_w)^{\flat_k}
= -h+\sum_{w\in\QQ}h|_k\gamma_w
= P_k(h)-2h,$$
and hence $\sigma_k(F)=h$. In particular, $F$ satisfies \textbf{(C3)}. It remains to verify \textbf{(C2)}. Let
$$\gamma_u=\begin{pmatrix}a&b\\c&d\end{pmatrix},
\qquad
\gamma_w=\begin{pmatrix}a'&b'\\c'&d'\end{pmatrix},$$
where $c,c'>0$. Then, modulo left multiplication by elements of $\langle T\rangle$, we have
\begin{equation}\label{eq:gamma-prod}
\gamma_w\gamma_u\equiv
\begin{cases}
\gamma_{(b+dw)/(a+cw)},&w>-a/c,\\
-\gamma_{(b+dw)/(a+cw)},&w<-a/c,\\
-I,&w=-a/c.
\end{cases}
\end{equation}
Moreover, since
$$\frac{b+dw}{a+cw}-u = -\frac1{c(a+cw)},$$
the change of variables $t=(b+dw)/(a+cw)$ maps $w<-a/c$ to $t>u$ and $w>-a/c$ to $t<u$. Thus, taking into account the signs arising from \cref{eq:gamma-prod} and the exceptional term $w=-a/c$, we obtain
$$P_k(h)|_k\gamma_u = \sum_{t\le u}h|_k\gamma_t - \sum_{t>u}h|_k\gamma_t-h.$$
Consequently,
$$F|_k\gamma_u
=\frac12P_k(h)|_k\gamma_u
=-F+\sum_{t\le u}h|_k\gamma_t
=-F+\sum_{t\le u}\sigma_k(F)|_k\gamma_t.$$
Thus, $F$ satisfies \textbf{(C2)}. Therefore, $F\in\EY_k$, proving surjectivity.
\end{proof}

\subsection{The even weight case}

\begin{lem}\label{lem:const-real}
Let $f\in\EY_k$ and $\alpha\ceq\lim_{\tau\to i\infty}f(\tau)$. Then, we have $\alpha\in\RR$.
\end{lem}

\begin{proof}
By \cref{lem:strong-cond}, we have $\sigma_k(f)\in \W_{k-3}$. Then, by \cref{eq:ord-3}, we have
$$P_k(\sigma_k(f))(\tau) - \sigma_k(f)(\tau)
= O\left(\sum_{c\geq 1,\,d\in\ZZ}\frac{1}{|c\tau+d|^3}\right) = O(\Img(\tau)^{-2})$$
as $\Img(\tau)\to\infty$. Hence, by \cref{lem:poincare},
$$(1-(-1)^k)\alpha = \lim_{\tau\to i\infty}P_k(\sigma_k(f))(\tau)
= \lim_{\tau\to i\infty}\sigma_k(f)(\tau) = \alpha-(-1)^k\overline{\alpha}.$$
It follows that $\alpha=\overline{\alpha}$, and hence $\alpha\in\RR$.
\end{proof}

\begin{lem}\label{lem:poincare-cusp}
Let $k$ be an even integer and $f\in\W_{k-3}^-$. Then, we have $P_k(f)\in iS_k$.
\end{lem}

\begin{proof}
By \cref{eq:ord-3}, the series defining $P_k(f)$ converges absolutely and locally uniformly on $\HH$, and hence defines a holomorphic function. Since $f$ is $1$-periodic and $k$ is even, $P_k(f)$ satisfies the modular transformation law. Moreover, by \cref{eq:ord-3}, we have
$$P_k(f)(\tau)=f(\tau)+O(\Img(\tau)^{-2})$$
as $\Img(\tau)\to\infty$, uniformly in $\operatorname{Re}(\tau)$. Thus, $P_k(f)$ is holomorphic at $i\infty$. Since $k$ is even, we have
$$(f|_k\gamma_w)^{\flat_k} = f^{\flat_k}|_k\gamma_{-w}.$$
Hence, using $f^{\flat_k}=-f$, we obtain
$$P_k(f)^{\flat_k} = P_k(f^{\flat_k}) = -P_k(f).$$
Thus, all Fourier coefficients of $P_k(f)$ are purely imaginary, so $P_k(f)\in iM_k$. Let $\alpha$ be the constant term of $f$. By the definition of $\W_{k-3}$, \cref{lem:const-real}, and the fact that $G(2)$ has a real constant term, we have $\alpha\in\RR$. Moreover,
$$\lim_{\tau\to i\infty}P_k(f)(\tau) = \lim_{\tau\to i\infty}f(\tau) = \alpha\in\RR.$$
Since $P_k(f)\in iM_k$, this implies that $\alpha\in\RR\cap i\RR$, and hence $\alpha=0$. Therefore, $P_k(f)\in iS_k$.
\end{proof}

\begin{thm}\label{thm:sigma-img-cusp}
For every even integer $k$, we have
$$P_k(\W_{k-3}^-)=iS_k.$$
More precisely, the forms $P_k(H_s)$ with $s$ odd and
$3\le s<(k-1)/2$ span $iS_k$ over $\RR$, where
$$H_s(\tau)\ceq G(s)(\tau)-\zeta(s)
= \frac{(-2\pi i)^s}{(s-1)!}\sum_{n\ge1}\sigma_{s-1}(n)q^n\qquad(q=e^{2\pi i\tau}).$$
\end{thm}

\begin{proof}
If $S_k=0$, the assertion follows from \cref{lem:poincare-cusp}. Thus, we may assume $S_k\ne0$, which implies $k\ge12$. Put
$$\mathcal I_k\ceq \{s\in\ZZ\mid s\text{ is odd},\ 3\le s<(k-1)/2\}.$$
For each $s\in\mathcal I_k$, we have $H_s\in\W_{k-3}^-$, and \cref{lem:poincare-cusp} gives $P_k(H_s)\in iS_k$. Let $d=\dim S_k$, and let $g_1,\ldots,g_d$ be an orthogonal basis of $S_k$ consisting of normalized Hecke eigenforms. Write
$$g_j(\tau)=\sum_{n\ge1}a_j(n)q^n.$$
By unfolding the Petersson inner product, we obtain
\begin{align*}
\langle P_k(H_s),g_j\rangle
&=\frac{(-2\pi i)^s}{(s-1)!}
\frac{\Gamma(k-1)}{(4\pi)^{k-1}}
\sum_{n\ge1}
\frac{\sigma_{s-1}(n)a_j(n)}{n^{k-1}}
=C_{k,s}
\frac{L(g_j,k-1)L(g_j,k-s)}{\zeta(k-s)},
\end{align*}
where
$$C_{k,s}\ceq \frac{(-2\pi i)^s}{(s-1)!}\frac{\Gamma(k-1)}{(4\pi)^{k-1}} \ne  0.$$
Here, the second equality follows from the Hecke relations. The condition $s<(k-1)/2$ ensures the absolute convergence of the relevant Dirichlet series. Moreover, $L(g_j,k-1)\ne0$ for every $j$, since the Euler product converges absolutely and is nonzero in the half-plane $\operatorname{Re}(z)>(k+1)/2$. Thus, the matrix
$$\bigl(\langle P_k(H_s),g_j\rangle
\bigr)_{j,s\in\mathcal I_k}$$
has the same rank as
$$\bigl(L(g_j,k-s)\bigr)_{j,s\in\mathcal I_k}.$$
We claim that this common rank is $d$. To prove this, it suffices to show that the $\RR$-linear map
$$S_k\longrightarrow\RR^{\mathcal I_k},
\qquad
g\longmapsto \bigl(L(g,k-s)\bigr)_{s\in\mathcal I_k}$$
is injective. Indeed, suppose that $g\in S_k$ satisfies $L(g,k-s)=0$ for every $s\in\mathcal I_k$. 
Here, we use the Eichler--Shimura theory of period polynomials (see \cite[Section 1]{KZ84} and \cite[Section 2]{Fuk07}).
Put $\til{k}=k-2$ and define the periods
$$r_n(g)\ceq\int_0^{i\infty}g(\tau)\tau^n\,d\tau = \frac{n!}{(-2\pi i)^{n+1}}L(g,n+1) \qquad(0\le n\le \til{k}).$$
By the functional equation of $L(g,z)$, the periods satisfy
$$r_n(g)+(-1)^nr_{\til{k}-n}(g)=0 \qquad(0\leq n\leq \til{k}).$$
The assumed vanishing of $L(g,k-s)$, together with these relations, implies
$$r_n(g)=0 \qquad(2\le n\le \til{k}-2,\ n\text{ even}).$$
Consequently, the even period polynomial
$$r_g^+(X)\ceq \sum_{\substack{0\le n\le \til{k}\\n\text{ even}}}\binom{\til{k}}{n}r_n(g)X^{\til{k}-n}$$
is a scalar multiple of $X^{\til{k}}-1$. By the Eichler--Shimura theory of even period polynomials, the map sending a cusp form to the class of its even period polynomial modulo $\mathbb C(X^{\til{k}}-1)$ is injective. Hence $g=0$, proving the claim. It follows that
$$\operatorname{rank}\bigl(L(g_j,k-s)\bigr)_{j,s\in\mathcal I_k} = d.$$
Therefore,
$$\operatorname{rank}\bigl(\langle P_k(H_s),g_j\rangle\bigr)_{j,s\in\mathcal I_k} = d.$$
Since $P_k(H_s)\in iS_k$ and the Petersson inner product is nondegenerate, the forms $P_k(H_s)$, $s\in\mathcal I_k$, span $iS_k$ over $\RR$. Therefore, we obtain
$$P_k(\W_{k-3}^-)=iS_k.$$
\end{proof}

\begin{lem}\label{lem:lift-kernel}
Let $k$ be an even integer and let
$$h\in\ker(P_k:\W_{k-3}^-\to iS_k).$$
Then there exists a holomorphic function $F_0$ on $\HH$, which is $1$-periodic and holomorphic at $i\infty$, such that
$$F_0|_k\gamma_u-F_0 = \sum_{w\le u}h|_k\gamma_w.$$
\end{lem}

\begin{proof}
For $\delta=\begin{pmatrix}a&b\\c&d\end{pmatrix}\in\SL_2(\ZZ)$, define
$$\rho_h(\delta)\ceq
\begin{cases}
\displaystyle\sum_{w\le d/c}h|_k\gamma_w,
&c\ne0,\\[2mm]
0,&c=0.
\end{cases}$$
We first show that $\rho_h$ satisfies the identity
\begin{equation}\label{eq:rho-cocycle}
\rho_h(\gamma\delta) = \rho_h(\gamma)|_k\delta+\rho_h(\delta)
\qquad(\gamma,\delta\in\SL_2(\ZZ)).
\end{equation}
To prove this, put $\gamma_\infty\ceq I$ and consider the M\"obius transformation
$$\psi_\delta(w)\ceq\frac{b+dw}{a+cw}$$
on $\PP^1(\QQ)$. By \cref{eq:gamma-prod}, right multiplication by $\delta$ induces the correspondence
$$\gamma_w\delta\equiv\gamma_{\psi_\delta(w)}$$
modulo left multiplication by elements of $\pm\langle T\rangle$, where the formula is understood on $\PP^1(\QQ)$. Moreover,
$$\psi_\delta'(w)=\frac1{(a+cw)^2}>0$$
wherever the derivative is defined. We regard $\PP^1(\RR)$ as an oriented circle, with the orientation induced by the increasing direction on $\RR$. For $p,q\in\PP^1(\RR)$, let $(p,q]$ denote the positively oriented arc from $p$ to $q$, with the initial point excluded and the terminal point included. An arc with identical endpoints is understood to be empty. For $\gamma\in\SL_2(\ZZ)$, let $x\in\PP^1(\QQ)$ be determined by $\gamma\equiv\gamma_x$ modulo left multiplication by $\pm\langle T\rangle$. Then
$$\rho_h(\gamma) = \sum_{w\in(\infty,x]\cap\PP^1(\QQ)}h|_k\gamma_w.$$
Since $\psi_\delta$ preserves the orientation of $\PP^1(\RR)$, right multiplication by $\delta$ maps the cosets indexed by $(\infty,x]$ bijectively onto those indexed by $(\psi_\delta(\infty),\psi_\delta(x)]$. Since $k$ is even, the signs of the coset representatives do not affect the slash action. Consequently,
$$\rho_h(\gamma)|_k\delta = \sum_{w\in(\psi_\delta(\infty),\psi_\delta(x)]\cap\PP^1(\QQ)}h|_k\gamma_w.$$
On the other hand, since
$$\sum_{w\in\PP^1(\QQ)}h|_k\gamma_w = P_k(h) = 0,$$
the sum over an oriented arc can be expressed as the difference of the corresponding sums starting at $\infty$. Indeed, any discrepancy is a multiple of the sum over the entire $\PP^1(\QQ)$, which vanishes. Therefore,
$$\rho_h(\gamma)|_k\delta
=\sum_{w\in(\infty,\psi_\delta(x)]\cap\PP^1(\QQ)}h|_k\gamma_w
- \sum_{w\in(\infty,\psi_\delta(\infty)]\cap\PP^1(\QQ)}h|_k\gamma_w
=\rho_h(\gamma\delta)-\rho_h(\delta).$$
This proves \cref{eq:rho-cocycle}. Next, consider the classical Poincar\'e series
$$Q\ceq \sum_{w\in\PP^1(\QQ)}q|_{12}\gamma_w
\qquad (q=e^{2\pi i\tau}).$$
Then $Q$ is a nonzero cusp form of weight $12$. Since $S_{12}=\CC\Delta$, we have $Q=C\Delta$ for some $C\ne0$. In particular, $Q$ has no zeros on $\HH$. Define
$$F_0\ceq-\frac1Q\sum_{w\in\PP^1(\QQ)}\rho_h(\gamma_w)(q|_{12}\gamma_w).$$
Since $\rho_h$ is invariant under left multiplication by elements of $\Gamma_\infty$, the sum is well-defined. By the uniform bounds above, it converges absolutely and locally uniformly on $\HH$. Hence, $F_0$ is holomorphic on $\HH$. For every $\delta\in\SL_2(\ZZ)$, \cref{eq:rho-cocycle} gives
\begin{align*}
F_0|_k\delta
&=-\frac{1}{Q}\sum_{w\in\PP^1(\QQ)}
(\rho_h(\gamma_w)|_k\delta)(q|_{12}\gamma_w\delta)
\\
&=-\frac{1}{Q}\sum_{w\in\PP^1(\QQ)}
\bigl(\rho_h(\gamma_w\delta)-\rho_h(\delta)\bigr)(q|_{12}\gamma_w\delta)
\\
&=F_0+\rho_h(\delta).
\end{align*}
Taking $\delta=T$ and $\delta=\gamma_u$ gives
$$F_0|_kT=F_0, \qquad
F_0|_k\gamma_u-F_0 = \sum_{w\le u}h|_k\gamma_w.$$
It remains to show that $F_0$ is holomorphic at $i\infty$. By \cref{eq:ord-3}, we have $\rho_h(\gamma_w)(\tau)=O(\Img(\tau)^{-2})$ uniformly in $w\in\QQ$ as $\Img(\tau)\to\infty$. Since $\rho_h(\gamma)=0$ when the lower-left entry of $\gamma$ is zero, we obtain
$$|Q(\tau)F_0(\tau)|
= O\left(\Img(\tau)^{-2}\sum_{c\ge1}\sum_{d\in\ZZ}|c\tau+d|^{-12}\right)
= O(\Img(\tau)^{-13})$$
as $\Img(\tau)\to\infty$. Since $QF_0$ is holomorphic and $1$-periodic, this estimate implies that its Fourier expansion has no constant or negative-power terms. Hence,
$$QF_0=O(q)\qquad(q=e^{2\pi i\tau}\to0).$$
On the other hand,
$$Q=C\Delta=Cq+O(q^2),\qquad C\ne0.$$
Therefore,
$$F_0=\frac{QF_0}{Q}=O(1)\qquad(q\to0).$$
Thus, $F_0$ is holomorphic at $i\infty$, completing the proof.
\end{proof}

\begin{thm}\label{thm:even-dim}
For every even integer $k$, we have
$$\sigma_k(\EY_k) = \ker(P_k:\W_{k-3}^-\to iS_k).$$
Consequently, by \cref{thm:sigma-img-cusp},
$$\dim\sigma_k(\EY_k) = \dim\W_{k-3}^- - \dim S_k.$$
\end{thm}

\begin{proof}
By \cref{lem:poincare}, we have $\sigma_k(\EY_k)\subseteq\ker P_k$. Conversely, let $h\in\W_{k-3}^-$ satisfy $P_k(h)=0$. By \cref{lem:lift-kernel}, there exists a holomorphic, $1$-periodic function $F_0$, holomorphic at $i\infty$, such that
$$F_0|_k\gamma_u-F_0 = \sum_{w\le u}h|_k\gamma_w\qquad(u\in\QQ).$$
Put $g\ceq \sigma_k(F_0)-h$. Since $F_0$ and $h$ are $1$-periodic, so is $g$.
We claim that $g|_kS=g$. Indeed, since $k$ is even and $h^{\flat_k}=-h$, we have
$$(h|_k\gamma_w)^{\flat_k} = h^{\flat_k}|_k\gamma_{-w} = -h|_k\gamma_{-w}.$$
Using $P_k(h)=0$, we obtain
\begin{align*}
(F_0|_kS-F_0)^{\flat_k}
= \left(\sum_{w\le0}h|_k\gamma_w\right)^{\flat_k}
= -\sum_{w\ge0}h|_k\gamma_w
= \sum_{w\le0}h|_k\gamma_w - h|_kS + h
= F_0|_kS-F_0 - h|_kS + h.
\end{align*}
Consequently,
$$g|_kS - g = (F_0|_kS-F_0) - (F_0|_kS-F_0)^{\flat_k} - (h|_kS-h) = 0.$$
Thus, $g$ is a modular form of weight $k$.
Moreover, $g^{\flat_k}=-g$, and hence $\sigma_k(g)=2g$. Define
$$F \ceq F_0-\frac{1}{2}g=\frac{1}{2}(F_0+F_0^{\flat_k}+h).$$
Then
$$\sigma_k(F)=\sigma_k(F_0)-g=h.$$
Since $g$ is a modular form, we also have
$$F|_k\gamma_u = F+\sum_{w\le u}h|_k\gamma_w
= F+\sum_{w\le u}\sigma_k(F)|_k\gamma_w.$$
The series converges absolutely by \cref{eq:ord-3}. Furthermore, $F$ is holomorphic, $1$-periodic, and holomorphic at $i\infty$, while $\sigma_k(F)=h\in\W_{k-3}^-$. Therefore $F\in\EY_k$, proving
$$\ker P_k\subseteq\sigma_k(\EY_k).$$
This completes the proof.
\end{proof}

\subsection{The dimension formula for $\EY_k$}

Combining the results for odd and even weights, we obtain the following exact sequence:
$$\begin{tikzcd}
0 \arrow[r] &
M_k \arrow[r, hook] &
\EY_k \arrow[r, "\sigma_k"] &
\W_{k-3}^- \arrow[r, "\til{P}_k"] &
iS_k \arrow[r] & 0,
\end{tikzcd}$$
where $\til{P}_k=P_k$ if $k$ is even and $\til{P}_k=0$ if $k$ is odd. Now, we define
$$\tEY_k\ceq\sum_{\substack{a,j\ge0\\2a+j=k}}G(2)^a\EY_j.$$
By \cref{thm:graded,thm:G2-indep}, this sum is direct, so
$$\tEY_k=\bigoplus_{\substack{a,j\ge0\\2a+j=k}}G(2)^a\EY_j.$$
Moreover, by \cref{thm:weight-graded,thm:inclusion}, we have $\E_k\subseteq\tEY_k$. Finally, by \cref{thm:product,thm:G2-indep}, the direct sum
$$\tEY\ceq\bigoplus_{k\ge0}\tEY_k$$
is a graded ring.

\begin{lem}\label{lem:dim-W-}
For every integer $k$, we have
$$\W_{k-3}^{-} \cong \bigoplus_{b\geq 0}(\W_{k-3-2b}/\W_{k-4-2b})
\cong \bigoplus_{b\geq 0}\tEY_{k-3-2b}.$$
\end{lem}

\begin{proof}
By \cref{thm:weight-graded,thm:inclusion,thm:G2-indep}, we have
$$\W_m=\bigoplus_{j=0}^m\tEY_j.$$
Consequently,
$$\W_m\cong\bigoplus_{j=0}^m\W_j/\W_{j-1},\qquad \W_j/\W_{j-1}\cong\tEY_j.$$
By \cref{lem:flat-mod-inv}, the involution $\flat_k$ acts on each quotient $\W_j/\W_{j-1}$ as multiplication by $(-1)^{k-j}$. Hence, we obtain
$$\W_{k-3}^-\cong\bigoplus_{\substack{0\le j\le k-3\\k-j\text{ odd}}}\W_j/\W_{j-1}.$$
The desired isomorphisms now follow by writing $j=k-3-2b$.
\end{proof}

\begin{lem}\label{lem:FY-dim}
For every integer $k$, we have
$$\dim \EY_{k} = \dim M_k + \sum_{n\geq 0} (n+1)\dim\EY_{k-3-2n} - \dim S_k.$$
\end{lem}

\begin{proof}
By the definition, we obtain
$$\dim(\W_{m}/\W_{m-1}) = \dim\tEY_m = \sum_{\substack{a,j\geq 0\\2a+j=m}}\dim\EY_{j}.$$
On the other hand, by \cref{thm:odd-dim,thm:even-dim}, we obtain
$$\dim\EY_k = \dim M_k + \dim\sigma_k(\EY_k) = \dim M_k + \dim\W_{k-3}^- - \dim S_k.$$
Therefore, using \cref{lem:dim-W-}, we have
\begin{align*}
\dim\EY_k &= \dim M_k + \sum_{\substack{a,b,j\geq 0\\2a+j=k-3-2b}}\dim\EY_{j} - \dim S_k
\\
&= \dim M_k + \sum_{n\geq 0} (n+1)\dim\EY_{k-3-2n} - \dim S_k
\qquad (\because n=a+b)
\end{align*}
This yields the desired equality.
\end{proof}

We are now ready to establish the dimension formula for $\EY_k$. Set
$$E(X)\ceq\frac{1}{1-X^2}, \qquad
O(X)\ceq\frac{X^3}{1-X^2},$$
$$M(X)\ceq\frac{1}{(1-X^4)(1-X^6)}, \qquad
S(X)\ceq\frac{X^{12}}{(1-X^4)(1-X^6)}.$$

\begin{thm}
The following coefficientwise inequality holds:
$$\sum_{k\ge 0}\dim\Ead_k X^k \preceq
\sum_{k\ge 0}\dim\EY_k X^k = \frac{M(X)-S(X)}{1-E(X)O(X)}.$$
\end{thm}

\begin{proof}
It is well known that
$$M(X)=\sum_{k\ge0}\dim M_k\,X^k
\quad\text{and}\quad
S(X)=\sum_{k\ge0}\dim S_k\,X^k.$$
By \cref{lem:FY-dim} and the identity
$$\frac{X^3}{(1-X^2)^2} = \sum_{n\ge0}(n+1)X^{2n+3},$$
we obtain
$$\sum_{k\ge0}\dim\EY_k\,X^k = M(X)-S(X)+\frac{X^3}{(1-X^2)^2} \sum_{k\ge0}\dim\EY_k\,X^k.$$
Rearranging this identity yields
$$\sum_{k\ge0}\dim\EY_k\,X^k
= \frac{M(X)-S(X)}{1-\frac{X^3}{(1-X^2)^2}}
= \frac{M(X)-S(X)}{1-E(X)O(X)},$$
as desired.
\end{proof}

\begin{cor}
The following coefficientwise inequality holds:
$$\sum_{k\ge 0}\dim\E_k X^k \preceq
\sum_{k\ge 0}\dim\tEY_k X^k = \frac{M(X)-S(X)}{1-X^2-O(X)}.$$
\end{cor}

\begin{proof}
By \cref{lem:FY-dim}, we obtain
$$\sum_{k\ge 0}\dim\tEY_k X^k
= \sum_{k\ge 0}\left(\dim\bigoplus_{2a+j=k}G(2)^a\EY_j\right) X^k = E(X)\sum_{k\ge0}\dim\EY_k\,X^k
= \frac{M(X)-S(X)}{1-X^2-O(X)},$$
as desired.
\end{proof}

Our original motivation was to characterize $\Ead_k$ in terms of modular transformation laws. However, the following dimension conjecture and the table comparing the dimensions of $\EY_k$ with the conjectural dimensions of $\Ead_k$ suggest that this goal cannot be achieved with the present definition. For low weights, we can verify that the conjectural dimensions of $\Ead_k$ give upper bounds for their actual dimensions, using the relations among multiple Eisenstein series introduced by \cite{BKM26}.

\begin{conj}[\cite{Oko14} with \cite{BKu14,BK26}] The following coefficientwise equality holds:
$$\sum_{k\geq 0}\dim\Ead_k X^k = \frac{M(X)}{1-E(X)O(X)+2E(X)S(X)}.$$
\end{conj}

\begin{table}[htbp]
\centering
\renewcommand{\arraystretch}{1.3}
\[
\begin{array}{c|*{18}{c}}
\hline
\text{weight }k
& 3 & 4 & 5 & 6 & 7 & 8 & 9 & 10 & 11 & 12
& 13 & 14 & 15 & 16 & 17 & 18 & 19 & 20
\\ \hline
\dim\EY_k
& 1 & 1 & 2 & 2 & 4 & 5 & 8 & 12 & 17 & 27
& 38 & 59 & 86 & 129 & 193 & 285 & 429 & 634
\\
\text{(Conj.)}\dim\Ead_k
& 1 & 1 & 2 & 2 & 4 & 5 & 8 & 12 & 17 & 26
& 38 & 57 & 83 & 124 & 183 & 271 & 401 & 593
\\ \hline
\end{array}
\]
\caption{Dimensions of $\EY_k$ and conjectural dimensions of $\Ead_k$.}
\label{tab:dimensions}
\end{table}

\bibliographystyle{myamsalpha}
\bibliography{ref}

@article{BKu14,
  author        = {Bachmann, Henrik and K{\"u}hn, Ulf},
  title         = {A short note on a conjecture of Okounkov about a {$q$}-analogue of multiple zeta values},
  journal       = {arXiv preprint},
  year          = {2014},
  eprint        = {1407.6796},
  archivePrefix = {arXiv},
  url           = {https://arxiv.org/abs/1407.6796}
}

@article{BK26,
  author        = {Bachmann, Henrik and Kanno, Hayato},
  title         = {The {$\mathfrak{sl}_2$}-algebra structure of multiple Eisenstein series},
  journal       = {arXiv preprint},
  year          = {2026},
  eprint        = {2609.03777},
  archivePrefix = {arXiv},
  url           = {https://arxiv.org/abs/2609.03777}
}

@article{BKM26,
  author        = {Bachmann, Henrik and Kanno, Hayato and Maesaka, Takumi},
  title         = {Relations and Derivatives of Multiple Eisenstein Series},
  journal       = {arXiv preprint},
  year          = {2026},
  eprint        = {2602.08176},
  archivePrefix = {arXiv},
  primaryClass  = {math.NT},
  url           = {https://arxiv.org/abs/2602.08176}
}

@article{Fuk07,
  author        = {Fukuhara, Shinji},
  title         = {Explicit formulas for Hecke operators on cusp forms, Dedekind symbols and period polynomials},
  journal       = {J. Reine Angew. Math.},
  volume        = {607},
  pages         = {163--216},
  year          = {2007}
}

@incollection{GKZ06,
  author        = {Gangl, Herbert and Kaneko, Masanobu and Zagier, Don},
  title         = {Double zeta values and modular forms},
  booktitle     = {Automorphic Forms and Zeta Functions},
  editor        = {B{\"o}cherer, Siegfried and Ibukiyama, Tomoyoshi and Kaneko, Masanobu and Sato, Fumihiro},
  publisher     = {World Scientific},
  address       = {Hackensack, NJ},
  pages         = {71--106},
  year          = {2006}
}

@incollection{KZ84,
  author        = {Kohnen, Winfried and Zagier, Don},
  title         = {Modular Forms with Rational Periods},
  booktitle     = {Modular Forms (Durham, 1983)},
  editor        = {Rankin, R. A.},
  series        = {Ellis Horwood Series in Mathematics and Its Applications},
  publisher     = {Ellis Horwood},
  address       = {Chichester},
  pages         = {197--249},
  year          = {1984}
}

@article{Oko14,
  author        = {Okounkov, Andrei},
  title         = {Hilbert schemes and multiple {$q$}-zeta values},
  journal       = {Funct. Anal. Appl.},
  volume        = {48},
  number        = {2},
  pages         = {138--144},
  year          = {2014}
}

@article{Yu26,
  author        = {Yu, J.},
  title         = {Modularity criteria for multiple Eisenstein series},
  journal       = {arXiv preprint},
  year          = {2026},
  eprint        = {2610.21291},
  archivePrefix = {arXiv},
  url           = {https://arxiv.org/abs/2610.21291}
}

\end{document}